\documentclass[12pt,oneside]{article}
\usepackage{amsfonts,pstricks}
\usepackage{amssymb}
\usepackage{euscript}
\usepackage{eufrak}
\usepackage{mathrsfs}
\usepackage{latexsym}
\usepackage{amsmath,enumerate}
\usepackage{amsthm}

\newtheorem{thm}{Theorem}[section]

\newtheorem{cor}{Corollary}[section]
\newtheorem{rmk}{Remark}[section]
\newtheorem{defn}{Definition}[section]
\newtheorem{prop}{Proposition}[section]
\newtheorem{ex}{Example}[section]

\newcommand{\NI}{\noindent}
\newcommand{\eb}{\vrule width 0.22 true cm height 0.22 true cm depth 0pt}

\usepackage{graphicx}
\date{}
\usepackage[normalem]{ulem}
\usepackage{hyperref}

\begin{document}

	\title
	{\large \bf Towards a Unified Framework for Bioperations on $e^\star$-Open Sets in Topological Spaces}
\author{G. Saravanakumar$^1$\footnote{saravananguru2612@gmail.com}\hspace{.2cm}, M. Arun$^2$\footnote{mathsarun@gmail.com} \footnote{$^{1}$Department of Mathematics, Vel Tech Rangarajan Dr.Sagunthala R\&D Institute of Science and Technology (Deemed to be University), Avadi, Chennai-600062, India., $^{2}$Department of Mathematics, Hindusthan College of Engineering \& Technology, Valley Campus. Coimbatore-	641032,India .}}
	\maketitle

	\begin{abstract}
		\NI \indent
		In this paper, we introduce the concept of $e^\star_{[\gamma,\gamma']}$-open sets in topological spaces and examine their properties in detail. Additionally, we propose a new class of functions, termed $(e^\star_{[\gamma,\gamma']},\ e^\star_{[\beta,\beta']})$-continuous functions, and explore their fundamental characteristics.
		\end{abstract}


	\textbf{AMS (2000) subject classification:}54A05,54A10

\section{Introduction}
Njastad \cite{3} introduced the concept of $\alpha$-open sets in topological spaces and investigated their properties. Kasahara \cite{2} defined the concept of an operation on topological spaces and introduced $\alpha$-closed graphs of an operation. Ogata \cite{4} referred to the operation $\alpha$ as the $\gamma$ operation and introduced the notion of $\tau_{\gamma}$, which is the collection of all $\gamma$-open sets in a topological space $(X,\ \tau)$. H. Z. Ibrahim \cite{6} extended the framework by defining operations on $\alpha O(X,\ \tau)$ and introduced $\alpha_{\gamma}$-open sets in topological spaces, exploring some of their fundamental properties.

Building upon these results, we extend the study to $e^\star$-open sets by employing a similar operation-based framework. E. Ekici \cite{eki} introduced the concept of $e^\star$-open sets and explored their fundamental properties in topological spaces. In this paper, we build on Ekici’s work and further analyze $e^\star$-open sets by integrating them into an extended operational framework. Below, we present the definitions relevant to $e^\star$-open sets and related concepts.

\section{Preliminaries}

Throughout this paper, $(X,\ \tau)$ and $(Y,\ \sigma)$ represent non-empty topological spaces on which no separation axioms are assumed, unless otherwise mentioned. Let $(X,\ \tau)$ be a topological space and $A$ a subset of $X$. The closure of $A$ and the interior of $A$ are denoted by $Cl(A)$ and $Int(A)$, respectively. A subset $A$ of a space $X$ is called regular open \cite{st} (resp. $\delta$-open \cite{veli}, $e^\star$-open \cite{eki}) if $A = Int(Cl(A)) \quad (\text{resp. } A = Int_\delta(A),\ A \subseteq Cl(Int(Cl_\delta(A)))).$
Here, $Int_\delta(A)$ is defined as $Int_\delta(A) = \{ x \mid U \in \mathcal{U}(x),\ Int(Cl(U)) \subseteq A \},$
where $\mathcal{U}(x)$ denotes the set of all open neighborhoods of the point $x$ in $X$. Equivalently, a subset $A$ of a space $X$ is called regular closed \cite{st} (resp. $\delta$-closed \cite{veli}, $e^\star$-closed \cite{eki}) if $A = Cl(Int(A)) \quad (\text{resp. } A = Cl_\delta(A),\ A \supseteq Int(Cl(Int_\delta(A)))).$
Here, $Cl_\delta(A)$ is given by $
Cl_\delta(A) = \{ x \mid \forall U \in \mathcal{U}(x),\ (Int(Cl(U)) \cap A) \neq \emptyset \}.$ The intersection of all $e^\star$- closed sets containing $A$ is called the $e^\star$-closure of $A$ and is denoted by $e^\star Cl(A)$. The family of all $e^\star$-open (resp. $e^\star$-closed) sets in a topological space $(X,\ \tau)$ is denoted by $e^\star O(X,\ \tau)$ (resp. $e^\star C(X, \tau)$ An operation $\gamma$ \cite{2} on a topology $\tau$ is a mapping from $\tau$ into power set $P(X)$ of $X$ such that $V\subseteq V^{\gamma}$ for each $V\in\tau,$ where $V^{\gamma}$ denotes the value of $\gamma$ at $V$. A subset $A$ of $X$ with an operation $\gamma$ on $\tau$ is called $\gamma$-open \cite{4} if for each $x\in A$, there exists an open set $U$ of $X$ containing $x$ such that $U^{\gamma}\subseteq A$. Clearly $\tau_{\gamma}\subseteq\tau$. Complements of $\gamma$-open sets are called $\gamma$-closed. An operation $\gamma$ : $\alpha O(X,\ \tau)\rightarrow P(X)$ \cite{6} is a mapping satisfying the condition, $V\subseteq V^{\gamma}$ for each $V\in \alpha O(X,\ \tau)$ . We call the mapping $\gamma$ an operation on $\alpha O(X,\ \tau)$ . A subset $A$ of $X$ is called an $\alpha_{\gamma}$-open set \cite{6} if for each point $x\in A$, there exists an $\alpha$-open set $U$ of $X$ containing $x$ such that $U^{\gamma}\subseteq A.$ We denote the set of all $\alpha_{\gamma}$-open sets of $(X,\ \tau)$ by $\alpha O(X,\ \tau)_{\gamma}$.
 An operation $\gamma$ on $\alpha O(X,\ \tau)$ is said to be $\alpha$-regular \cite{6} if for every $\alpha$-open sets $U$ and $V$ containing $x\in X$, there exists an $\alpha$-open set $W$ of $x$ such that $W^{\gamma}\subseteq U^{\gamma}\cap V^{\gamma}.$ An operation $\gamma$ on $\alpha O(X,\ \tau)$ is said to be $\alpha$-open \cite{6} if for every $\alpha$-open set $U$ of $x\in X$, there, exists an $\alpha_{\gamma}$-open set $V$ of $X$ such that $x\in V$ and $V\subseteq U^{\gamma}.$ A non-empty subset \( A \) of \( (X, \tau) \) is said to be \( \alpha_{[\gamma,\gamma']} \)-open\cite{7a} if \( \forall x \in A \), \( \exists \) \( \alpha \)-open sets \( U \) and \( V \) of \( X \) containing \( x \), such that \( U^{\gamma} \cap V^{\gamma'} \subseteq A \). The set of all \( \alpha_{[\gamma,\gamma']} \)-open sets of \( (X, \tau) \) is denoted by \( \alpha O(X, \tau)_{[\gamma,\gamma']} \). If \( A_i \) is \( \alpha_{[\gamma,\gamma']} \)-open \( \forall i \in I \), then \( \bigcup_{i \in I} A_i \) is \( \alpha_{[\gamma,\gamma']} \)-open\cite{7a}.

\section{$e^\star_{[\gamma,\gamma']}$-Open Sets}

\defn\em A non-empty subset \( A \) of \( (X, \tau) \) is said to be \( e^\star_{[\gamma,\gamma']} \)-open if \( \forall x \in A \), \( \exists \) \( e^\star \)-open sets \( U \) and \( V \) of \( X \) containing \( x \), such that \( U^{\gamma} \cap V^{\gamma'} \subseteq A \). The set of all \( e^\star_{[\gamma,\gamma']} \)-open sets of \( (X, \tau) \) is denoted by \( e^\star O(X, \tau)_{[\gamma,\gamma']} \).

We suppose that the empty set is $e^\star_{[\gamma,\gamma']}$-open for any operations $\gamma$ and $\gamma'.$

\defn\em An operation $\gamma$ on $e^\star O(X, \tau)$ is said to be $e^\star$-regular if for every $e^\star$-open sets $U$ and $V$ containing $x\in X$, there exists an $e^\star$-open set $W$ of $x$ such that $W^{\gamma}\subseteq U^{\gamma}\cap V^{\gamma}.$

\prop\em\label{p3.1} If \( A_i \) is \( e^\star_{[\gamma,\gamma']} \)-open \( \forall i \in I \), then \( \bigcup_{i \in I} A_i \) is \( e^\star_{[\gamma,\gamma']} \)-open.

\noindent\textbf{Proof.} Let $x\displaystyle \in\bigcup_{i\in I}A_{i}$, then $x\in A_{i}$ for some $i\in I$. Since $A_{i}$ is an $e^\star_{[\gamma,\gamma']}$-open set, so there exist $e^\star$-open sets $U$ and $V$ of $(X, \tau)$ such that $ x\in U\cap V\subseteq U^{\gamma}\displaystyle \cap V^{\gamma'}\subseteq A_{i}\subseteq\bigcup_{i\in I}A_{i}.$ Therefore $\displaystyle \bigcup_{i\in I}A_{i}$ is an $e^\star_{[\gamma,\gamma']}$-open set of $(X, \tau)$ . \eb

If $A$ and $B$ are two $e^\star_{[\gamma,\gamma']}$-open sets in $(X, \tau)$ , then the following example shows that $A\cap B$ need not be $e^\star_{[\gamma,\gamma']}$-open.

\ex\em\label{e3.1} Let $X=\{w_{1}, w_{2}, w_{3}\}$  and  $\tau=\{\phi, X, \{w_2\}, \{w_3\}, \{w_1,\ w_2\},$ $ \{w_2, w_3\}\}$ be a topology on $X$. For each $A\in e^\star O(X, \tau)$, we define two operations $\gamma$ and $\gamma',$ by

$A^{\gamma}=\left\{\begin{array}{ll}
Cl(A) & \mathrm{i}\mathrm{f}\ w_{3}\in A,\\
X & \mathrm{i}\mathrm{f}\ w_{3}\not\in A,
\end{array}\right.$ and $A^{\gamma'}=\left\{\begin{array}{ll}
A & \mathrm{i}\mathrm{f}\ A\neq\{w_{2}\},\\
X & \mathrm{i}\mathrm{f}\ A=\{w_{2}\}.
\end{array}\right.$

Then, it is obvious that the sets $\{w_{1}, w_{2}\}$ and $\{w_{2},w_{3}\}$ are $e^\star_{[\gamma,\gamma']}$-open, however their intersection $\{w_{2}\}$ is not $e^\star_{[\gamma,\gamma']}$-open.

\prop\em\label{p3.2} Let $\gamma$ and $\gamma'$ be $e^\star$-regular operations. If $A$ and $B$ are $e^\star_{[\gamma,\gamma']}$-open, then $A\cap B$ is $e^\star_{[\gamma,\gamma']}$-open.

\noindent\textbf{Proof.} Let $x\in A\cap B$, then $x\in A$ and $x\in B$. Since $A$ and $B$ are $e^\star_{[\gamma,\gamma']^{-}}$ open sets, there exist $e^\star$-open sets $U, V, W$ and $S$ such that $ x\in U\cap V\subseteq U^{\gamma}\cap V^{\gamma'}\subseteq A$ and $x\in W\cap S\subseteq W^{\gamma}\cap S^{\gamma'}\subseteq B$. Since $\gamma$ and $\gamma'$ are $e^\star$-regular operations, then there exist an $e^\star$-open sets $K$ and $L$ containing $x$ such that $K^{\gamma}\cap L^{\gamma'}\subseteq(U^{\gamma}\cap W^{\gamma})\cap(V^{\gamma'}\cap S^{\gamma'})=(U^{\gamma}\cap V^{\gamma'})\cap(W^{\gamma}\cap S^{\gamma'})\subseteq A\cap B$. This implies that $A\cap B$ is $e^\star_{[\gamma,\gamma']}$-open set. \eb

\rmk\em By the above propositon, if $\gamma$ and $\gamma'$ are $e^\star$-regular operations, then $e^\star O(X, \tau)_{[\gamma,\gamma']}$ forms a topology on $X.$

\prop\em The set $A$ is $e^\star_{[\gamma,\gamma']}$-open in $X$ if and only if for each $x \in A$, there exists an $e^\star_{[\gamma,\gamma']}$-open set $B$ such that $x \in B \subseteq A.$

\noindent\textbf{Proof.} Suppose that $A$ is an $e^\star_{[\gamma,\gamma']}$-open set in the space $X$. Then for each $x\in A$, put $B=A$ which is an $e^\star_{[\gamma,\gamma']}$-open set such that $x\in B\subseteq A.$

Conversely, suppose that for each $x\in A$, there exists an $e^\star_{[\gamma,\gamma']}$-open set $B$ such that $x\in B\subseteq A$. Thus $A=\displaystyle \bigcup_{x\in A}B_{x}$, where $B_{x}\in e^\star O(X,\ \tau)_{[\gamma,\gamma']}$. Therefore, $A$ is an $e^\star_{[\gamma,\gamma']}$-open set. \eb

\prop\em If $A$ is $e^\star_{[\gamma,\gamma']}$-open, then $A$ is $e^\star$-open.

The converse of the above proposition need not be true in general as it is shown below.

\ex\em Let $(X, \tau), \gamma$ and $\gamma'$ be the same space and the same operations as in Example \ref{e3.1} Then $\{w_{2}\}$ is $e^\star$-open but not $e^\star_{[\gamma,\gamma']}$-open.

\rmk\em\label{r3.3} A subset $A$ is an $e^\star_{[id,id']}$-open set of $(X, \tau)$ if and only if $A$ is $e^\star$-open in $(X, \tau)$ . The operation $id:e^\star O(X, \tau)\rightarrow P(X)$ is defined by $V^{id}=V$ for any set $V\in e^\star O(X, \tau)$ . This operation is called the identity operation on $e^\star O(X, \tau)$ . Therefore $e^\star O(X, \tau)_{[id,id']}=e^\star O(X, \tau)$ .

\rmk\em From Remark \ref{r3.3} we have $e^\star O(X, \tau)_{[id,id']}= e^\star O(X, \tau)=e^\star O(X, \tau)_{id}=e^\star O(X, \tau)_{id'}.$

\rmk\em The following example shows that the concept of $e^\star_{[\gamma,\gamma']^{-}}$ open and open are independent.

\ex\em\label{33} Let $X=\{u_{1}, u_{2}, u_{3}\}$ and $\tau=\{\phi, X, \{u_{1}\}\},$ be a topology on $X$. For each $A\in e^\star O(X,\ \tau)$ we define two operations $\gamma$ and $\gamma'$, respectively, by $A^{\gamma}=X$ and $A^{\gamma'}=\left\{\begin{array}{ll}
A & \mathrm{i}\mathrm{f}\ A=\{u_{1}, u_{2}\},\\
X & \mathrm{i}\mathrm{f}\ A\neq\{u_{1}, u_{2}\}.
\end{array}\right.$

Then, $e^\star_{[\gamma,\gamma']}$-open sets are $\phi, X$, and $\{u_{1}, u_{2}\}.$

\prop\em If $A$ is $[\gamma, \gamma']$-open, then $A$ is $e^\star_{[\gamma,\gamma']}$-open.

The converse of the above proposition need not be true in general as it is shown below.

\ex\em Let $X=\{v_{1}, v_{2}, v_{3}\}$ and $\tau=\{\phi, X, \{v_{2}\}\}$ be a topology on { X}. For each $A\in e^\star O(X, \tau)$, we define two operations $\gamma$ and $\gamma'$, respectively, by $A^{\gamma}=A^{\gamma'}=A$. Then $\{v_{1}, v_{2}\}$ is $e^\star_{[\gamma,\gamma']}$-open but not $[\gamma, \gamma']$-open.

\prop\em\label{p3.6} If $A$ is $e^\star_{\gamma}$-open and $B$ is $e^\star_{\gamma'}$-open, then $A \cap B$ is $e^\star_{[\gamma,\gamma']}$-open.

\prop\em\label{p3.7} If $A$ is $e^\star_{\gamma}$-open, then $A$ is $e^\star_{[\gamma,\gamma']}$-open for any operation $\gamma'$.

The converse of Proposition \ref{p3.7} need not be true in general as it is shown below.

\ex\em Let $X=\{t_{1}, t_{2}, t_{3}\}$ and $\tau$ be a discrete topology on $X$. For each $A\in e^\star O(X, \tau)$ we define two operations $\gamma$ and $\gamma'$, respectively, by $A^{\gamma}=A^{\gamma'}=\left\{\begin{array}{ll}
A & \mathrm{i}\mathrm{f}\ A=\{t_{1}, t_{2}\}~\text{or}~\{t_{2}, t_{3}\},\\
X & \text{otherwise}.
\end{array}\right.$

Then $\{t_{2}\}$ is $e^\star_{[\gamma,\gamma']}$-open but not $e^\star_{\gamma}$-open.

\prop\em Let $X: e^\star O(X, \tau) \rightarrow P(X)$ be an operation defined by $U^{X} = X$ for every $U \in e^\star O(X, \tau)$. Then $A$ is $e^\star_{\gamma}$-open if and only if $A$ is $e^\star_{[\gamma, \gamma']}$-open.

\defn\em A topological space $(X,\ \tau)$ is said to be $e^\star_{\gamma}$-regular if for each $x\in X$ and for each $e^\star$-open set $V$ in $X$ containing $x$, there exists an $e^\star$-open set $U$ in $X$ containing $x$ such that $U^{\gamma}\subseteq V.$

In the following proposition, we give a condition under which the family of $e^\star$-open sets is identical to the family of $e^\star_{\gamma}$-open sets.

\prop\em A topological space $(X, \tau)$ with an operation $\gamma$ on $e^\star O(X, \tau)$ is $e^\star_{\gamma}$-regular if and only if $e^\star O(X, \tau) = e^\star O(X, \tau)_{\gamma}.$

\rmk\em If a topological space $(X, \tau)$ is $e^\star_{\gamma}$-regular, then $\tau\subseteq e^\star O(X, \tau)_{\gamma}.$

\defn\em A topological space $(X, \tau)$ is said to be $e^\star_{[\gamma,\gamma']}$-regular if for each point $x$ in $X$ and every $e^\star$-open set $U$ containing $x$ there exist $e^\star$-open sets $W$ and $S$ of $x$ such that $W^{\gamma}\cap S^{\gamma'}\subseteq U.$

\prop\em\label{p3.10} A topological space $(X, \tau)$ with operations $\gamma$ and $\gamma'$ on $e^\star O(X, \tau)$ is $e^\star_{[\gamma, \gamma']}$-regular if and only if $e^\star O(X, \tau) = e^\star O(X, \tau)_{[\gamma, \gamma']}.$

\noindent\textbf{Proof.} Let $(X, \tau)$ be $e^\star_{[\gamma,\gamma']}$-regular and $A\in e^\star O(X, \tau)$. Since $(X, \tau)$ is $e^\star_{[\gamma,\gamma']}$-regular, then for each $x\in A$, there exist $e^\star$-open sets $W$ and $S$ of $x$ such that $W^{\gamma}\cap S^{\gamma'}\subseteq A$. This implies that $A\in e^\star O(X, \tau)_{[\gamma,\gamma']}$. But we have $e^\star O(X, \tau)_{[\gamma,\gamma']}\subseteq e^\star O(X, \tau)$. Therefore $e^\star O(X, \tau)=e^\star O(X, \tau)_{[\gamma,\gamma']}.$

Conversely, let $e^\star O(X, \tau)=e^\star O(X, \tau)_{[\gamma,\gamma']}, x\in X$ and $V$ be $e^\star$-open of $x$. Then by assumption $V$ is $e^\star_{[\gamma,\gamma']}$-open set. This implies that there exist $e^\star$-open sets $W$ and $S$ of $x$ such that $W^{\gamma}\cap S^{\gamma'}\subseteq V$. Therefore $(X, \tau)$ is $e^\star_{[\gamma,\gamma']}$-regular.\eb

\prop\em For $e^\star_{\gamma}$-regularity, $e^\star_{\gamma'}$-regularity, and $e^\star_{[\gamma, \gamma']}$-regularity of a space $(X, \tau)$, the following properties hold.

\begin{enumerate}[(1)]
\item  $(X, \tau)$ is $e^\star_{[\gamma, X]}$-regular if and only if it is $e^\star_{\gamma}$-regular.
\item If $(X, \tau)$ is $e^\star_{\gamma}$-regular and $e^\star_{\gamma'}$-regular, then it is $e^\star_{[\gamma, \gamma']}$-regular.
\end{enumerate}

\defn\em\label{d34} A subset $F$ of $(X, \tau)$ is said to be $e^\star_{[\gamma,\gamma']}$-closed if its complement $X\backslash F$ is $e^\star_{[\gamma,\gamma']}$-open.

We denote the set of all $e^\star_{[\gamma,\gamma']}$-closed sets of $(X, \tau)$ by $e^\star C(X, \tau)_{[\gamma,\gamma']}.$

\defn\em\label{d35} Let $A$ be a subset of a topological space $(X, \tau)$ . The intersection of all $e^\star_{[\gamma,\gamma']}$-closed sets containing $A$ is called the $e^\star_{[\gamma,\gamma']}$-closure of $A$ and denoted by $e^\star_{[\gamma,\gamma']}$-$Cl(A)$ .

\prop\em\label{p3.12} For a point $x \in X$, $x \in e^\star_{[\gamma, \gamma']}$-$Cl(A)$ if and only if $V \cap A \neq \emptyset$ for every $e^\star_{[\gamma, \gamma']}$-open set $V$ containing $x.$

\prop\em\label{p3.13} Let $A$ and $B$ be subsets of $(X, \tau)$. Then the following hold:
\begin{enumerate}[(1)]
	\item  $A\subseteq e^\star_{[\gamma,\gamma']}$-$Cl(A)$.
	\item If $A\subseteq B$, then $e^\star_{[\gamma,\gamma']}$-$Cl(A)\subseteq e^\star_{[\gamma,\gamma']}$-$Cl(B)$.
	\item  $A\in e^\star C(X, \tau)_{[\gamma,\gamma']}$ if and only if $e^\star_{[\gamma,\gamma']}$-$Cl(A)=A.$
	\item  $e^\star_{[\gamma,\gamma']}$-$Cl(A)\in e^\star C(X, \tau)_{[\gamma,\gamma']}.$
	\item  $e^\star_{[\gamma,\gamma']}$-$Cl(A\cap B)\subseteq e^\star_{[\gamma,\gamma']}$-$Cl(A)\cap e^\star_{[\gamma,\gamma']}$-$Cl(B).$
	\item  If $\gamma$ and $\gamma'$ are $e^\star$-regular, then $e^\star_{[\gamma,\gamma']}$-$Cl(A\cup B)=e^\star_{[\gamma,\gamma']}$-$Cl(A)\cup e^\star_{[\gamma,\gamma']}$-$Cl(B).$
\end{enumerate}

We introduce the following definition of $e^\star Cl_{[\gamma,\gamma']}(A)$ .

\defn\em\label{d36} For a subset $A$ of $(X, \tau)$ , we define $e^\star Cl_{[\gamma,\gamma']}(A)$ as follows: $e^\star Cl_{[\gamma,\gamma']}(A)=\{x\in X$ : $(U^{\gamma}\cap W^{\gamma'})\cap A\neq\phi$ holds for every $e^\star$-open sets $U$ and $W$ containing $x$\}.

\rmk\em In Definitions \ref{d34},\ref{d35} and \ref{d36}, put $\gamma'=X$. Then, for any subset $A$ of $X$, the following hold:

\begin{enumerate}[(1)]
\item  $e^\star_{[\gamma,X]}$-$Cl$ $(A)=e^\star_{\gamma}$-$Cl(A)$ .
\item  $e^\star C(X, \tau)_{[\gamma,X]}=$\{ $F:F$ is $e^\star_{\gamma}$-closed\}.
\item  $e^\star Cl_{[\gamma,X]}(A)=e^\star Cl_{\gamma}(A)$ .
\end{enumerate}

\prop\em\label{p3.4} For a subset $A$ of $(X, \tau)$, we have $A \subseteq e^\star Cl(A) \subseteq e^\star Cl_{[\gamma, \gamma']}(A) \subseteq e^\star_{[\gamma, \gamma']}$-$Cl(A)$.

\thm\em\label{t3.1} Let $A$ be a subset of a topological space $(X, \tau)$, the following properties are equivalent:
\begin{enumerate}[(1)]
	\item  $A\in e^\star O(X, \tau)_{[\gamma,\gamma']}.$
	\item  $e^\star Cl_{[\gamma,\gamma']}(X\backslash A)=X\backslash A.$
	\item  $e^\star_{[\gamma,\gamma']}$-$Cl(X\backslash A)=X\backslash A.$
	\item  $X\backslash A\in e^\star C(X, \tau)_{[\gamma,\gamma']}.$
\end{enumerate}

\thm\em For a subset $A$ of $(X, \tau)$, the following properties hold:
\begin{enumerate}[(1)]
	\item If $(X, \tau)$ is $e^\star_{[\gamma,\gamma']}$-regular, then $e^\star Cl(A)=e^\star Cl_{[\gamma,\gamma']}(A)=e^\star_{[\gamma,\gamma']}$-$Cl(A)$.
	\item $e^\star Cl_{[\gamma,\gamma']}(A)$ is an $e^\star$-closed subset of $(X, \tau)$.
\end{enumerate}

\thm\em\label{t3.3} Let $\gamma$ and $\gamma'$ be $e^\star$-open operations and $A$ a subset of $(X, \tau)$. Then, the following statements hold:
\begin{enumerate}[(1)]
	\item  $e^\star Cl_{[\gamma, \gamma']}(A)=e^\star_{[\gamma, \gamma']}$-$Cl(A).$
	\item  $e^\star Cl_{[\gamma, \gamma']}(e^\star Cl_{[\gamma, \gamma']}(A))=e^\star Cl_{[\gamma, \gamma']}(A).$
\end{enumerate}

\noindent\textbf{Proof.} (1). By Proposition \ref{p3.4}, it suffices to prove that $e^\star_{[\gamma,\gamma']}$-$Cl(A)\subseteq e^\star Cl_{[\gamma,\gamma']}(A)$ . Let $x\in e^\star_{[\gamma,\gamma']}$-$Cl(A)$ and $W$ and $S$ be $e^\star$-open sets of $X$ containing $x$. By the $e^\star$-openness of $\gamma$ and $\gamma'$, there exist an $e^\star_{\gamma}$-open set $W'$ and an $e^\star_{\gamma'}$-open set $S'$ such that $x\in W'\subseteq W^{\gamma}$ and $x\in S'\subseteq S^{\gamma'}$. By Propositions \ref{p3.6} and \ref{p3.12}, $(S'\cap W')\cap A\neq\phi$ and hence $(S^{\gamma}\cap W^{\gamma'})\cap A\neq\phi.$ This implies that $x\in e^\star Cl_{[\gamma,\gamma']}(A)$ .

2. This follows immediately from (1) and Proposition \ref{p3.13} (3).\eb

\rmk\em The below example shows that the equalities of Theorem \ref{t3.3} are not true without the assumption that both operations are $e^\star$-open.

\ex\em Let $X=\{z_1, z_2, z_3\}$ and $\tau=\{\phi, X, \{z_1\}, \{z_2\}, \{z_1, z_2\}\}$ For each $A\in e^\star O(X,\ \tau)$ we define two operations $\gamma$ and $\gamma'$, respectively, by $A^{\gamma}= Cl(A)$ and $A^{\gamma'}=X$. The operation $\gamma$ is not $e^\star$-open. However, $e^\star Cl_{[\gamma,\gamma']}\{z_1\}= \{z_1,\ z_3\}\subseteq e^\star_{[\gamma,\gamma']}$-$Cl(\{z_1\})=X$ and $e^\star Cl_{[\gamma,\gamma']}(e^\star Cl_{[\gamma,\gamma']}(\{z_1\}))=X\neq\{z_1, z_3\}= e^\star Cl_{[\gamma,\gamma']}(\{z_1\})$ .

\thm\em\label{t3.4}Let $A$ and $B$ be subsets of a topological space $(X, \tau)$ and $\gamma, \gamma'$: $e^\star O(X, \tau) \rightarrow P(X)$ operations on $e^\star O(X, \tau)$. Then we have the following properties:
\begin{enumerate}[(1)]
	\item  $A \subseteq e^\star Cl_{[\gamma,\gamma']}(A)$.
	\item  $e^\star Cl_{[\gamma,\gamma']}(\emptyset) = \emptyset$ and $e^\star Cl_{[\gamma,\gamma']}(X) = X.$
	\item  $A \in e^\star C(X, \tau)_{[\gamma,\gamma']}$ if and only if $e^\star Cl_{[\gamma,\gamma']}(A) = A.$
	\item If $A \subseteq B$, then $e^\star Cl_{[\gamma,\gamma']}(A) \subseteq e^\star Cl_{[\gamma,\gamma']}(B).$
	\item  $e^\star Cl_{[\gamma,\gamma']}(A \cup B) \subseteq e^\star Cl_{\gamma}(A) \cup e^\star Cl_{\gamma'}(B).$
	\item  If $\gamma$ and $\gamma'$ are $e^\star$-regular, then $e^\star Cl_{[\gamma,\gamma']}(A \cup B) = e^\star Cl_{[\gamma,\gamma']}(A) \cup e^\star Cl_{[\gamma,\gamma']}(B).$
	\item  $e^\star Cl_{[\gamma,\gamma']}(A \cap B) \subseteq e^\star Cl_{[\gamma,\gamma']}(A) \cap e^\star Cl_{[\gamma,\gamma']}(B).$
\end{enumerate}

\textbf{Proof.} (1), (2) and (4). They are obtained from Definition \ref{d36}.

(3). Suppose that $A$ is $e^\star_{[\gamma,\gamma']}$-closed. Then $X\backslash A$ is $e^\star_{[\gamma,\gamma']}$-open in $(X,\ \tau)$ . We claim that $e^\star Cl_{[\gamma,\gamma']}(A)\subseteq A$. Let $x\not\in A$. There exist $e^\star$-open sets $U$ and $V$ of $(X,\ \tau)$ containing $x$ such that $U^{\gamma}\cap V^{\gamma'}\subseteq X\backslash A$, that is, $(U^{\gamma}\cap V^{\gamma'})\cap A=\phi$. Hence by Definition \ref{d36}, we have that $x\not\in e^\star Cl_{[\gamma,\gamma']}(A)$ and so $e^\star Cl_{[\gamma,\gamma']}(A)\subseteq A$. By (1), it is proved that $e^\star Cl_{[\gamma,\gamma']}(A)=A$. Conversely, suppose that $e^\star Cl_{[\gamma,\gamma']}(A)= A$. Let $x\in X\backslash A$. Since $x\not\in e^\star Cl_{[\gamma,\gamma']}(A)$ , there exist an $e^\star$-open sets $U$ and $V$ containing $x$ such that $(U^{\gamma}\cap V^{\gamma'})\cap A=\phi$, that is, $U^{\gamma}\cap V^{\gamma'}\subseteq X\backslash A.$ Therefore, $A$ is $e^\star_{[\gamma,\gamma']}$-closed.

(5), (7). They are obtained from (4).

(6). Let $x\not\in e^\star Cl_{[\gamma,\gamma']}(A)\cup e^\star Cl_{[\gamma,\gamma']}(B)$ . Then there exist $e^\star$-open sets $U, V, W$ and $S$ of $(X,\ \tau)$ containing $x$ such that $(U^{\gamma}\cap V^{\gamma'})\cap A=\phi$ and $(W^{\gamma}\cap S^{\gamma'})\cap B=\phi.$ Since $\gamma$ and $\gamma'$ are $e^\star$-regular, by definition of $e^\star$-regular, there exist $e^\star$-open sets $K$ and $L$ of $(X,\ \tau)$ containing $x$ such that $K^{\gamma}\subseteq U^{\gamma}\cap W^{\gamma}$ and $L^{\gamma'}\subseteq V^{\gamma'}\cap S^{\gamma'}$ Thus, we have $(k^{\gamma}\cap L^{\gamma'})\cap(A\cup B)\subseteq((U^{\gamma}\cap W^{\gamma})\cap(V^{\gamma'}\cap S^{\gamma})) \cap(A\cup B)=((U^{\gamma}\cap  V^{\gamma'})\cap(W^{\gamma}\cap S^{\gamma'}))\cap(A\cup B)=[((U^{\gamma}\cap V^{\gamma'})\cap(W^{\gamma}\cap S^{\gamma'}))\cap A]\cup[((U^{\gamma}\cap V^{\gamma'})\cap(W^{\gamma}\cap  S^{\gamma'}))\cap B]=\phi$, that is, $(K^{\gamma}\cap L^{\gamma'})\cap(A\cup B)=\phi$. Hence, $x\not\in e^\star Cl_{[\gamma,\gamma']}(A\cup B)$ . This shows that $e^\star Cl_{[\gamma,\gamma']}(A)\cup e^\star Cl_{[\gamma,\gamma']}(B)\supseteq e^\star Cl_{[\gamma,\gamma']}(A\cup B)$ . \eb

\rmk\em Example \ref{e3.1} shows that the inclusion of Theorem \ref{t3.4} (5) is a proper one in general. For a subset $\{w_1\}, e^\star Cl_{[\gamma,\gamma']}(\{w_1\})=\{a\}\subseteq e^\star Cl_{\gamma}(\{w_1\})\cup e^\star Cl_{\gamma'}(\{w_1\})=\{w_1,\ w_2 \}.$

We define the $e^\star_{[\gamma,\gamma']}$-interior of a subset $A$ of $(X,\ \tau)$ as follows:

\defn\em Let $A$ be a subset of a topological space $(X,\ \tau)$ . The union of all $e^\star_{[\gamma,\gamma']}$-open sets contained in $A$ is called the $e^\star_{[\gamma,\gamma']}$-interior of $A$ and is denoted by $e^\star_{[\gamma,\gamma']}$-$Int(A)$ .

\prop\em For any subsets $A, B$ of $X$, we have the following:
\begin{enumerate}[(1)]
	\item $e^\star_{[\gamma,\gamma']}$-Int(A) is an $e^\star_{[\gamma,\gamma']}$-open set in $X.$
	\item $A$ is $e^\star_{[\gamma,\gamma']}$-open if and only if $A = e^\star_{[\gamma,\gamma']}$-$Int(A).$
	\item $e^\star_{[\gamma,\gamma']}$-$Int(e^\star_{[\gamma,\gamma']}$-$Int(A)) = e^\star_{[\gamma,\gamma']}$-$Int(A).$
	\item $e^\star_{[\gamma,\gamma']}$-$Int(A) \subseteq A.$
	\item If $A \subseteq B$, then $e^\star_{[\gamma,\gamma']}$-$Int(A) \subseteq e^\star_{[\gamma,\gamma']}$-$Int(B)$.
	\item $e^\star_{[\gamma,\gamma']}$-$Int(A \cup B) \supseteq e^\star_{[\gamma,\gamma']}$-$Int(A) \cup e^\star_{[\gamma,\gamma']}$-$Int(B).$
	\item $e^\star_{[\gamma,\gamma']}$-$Int(A \cap B) \subseteq e^\star_{[\gamma,\gamma']}$-$Int(A) \cap e^\star_{[\gamma,\gamma']}$-$Int(B).$
\end{enumerate}

\textbf{Proof.} Obvious. \eb

\prop\em Let $A$ be any subset of a topological space $(X, \tau)$. Then the following statements are true:
\begin{enumerate}[(1)]
	\item $X\backslash e^\star_{[\gamma, \gamma']}$-$Int(A) = e^\star_{[\gamma, \gamma']}$-$Cl(X\backslash A).$
	\item $X\backslash e^\star_{[\gamma, \gamma']}$-$Cl(A) = e^\star_{[\gamma, \gamma']}$-$Int(X\backslash A).$
	\item $e^\star_{[\gamma, \gamma']}$-$Int(A) = X\backslash e^\star_{[\gamma, \gamma']}$-$Cl(X\backslash A).$
	\item $e^\star_{[\gamma, \gamma']}$-$Cl(A) = X\backslash e^\star_{[\gamma, \gamma']}$-$Int(X\backslash A).$
\end{enumerate}

\section{$(e^\star_{[\gamma,\gamma']},\ e^\star_{[\beta,\beta']})$-Functions}

Throughout this section, let $f$ : $(X,\ \tau)\rightarrow(Y,\ \sigma)$ be a function and $\gamma, \gamma'$ : $e^\star O(X,\ \tau)\rightarrow P(X)$ be operations on $e^\star O(X,\ \tau)$ and $\beta, \beta'$ : $e^\star O(Y,\ \sigma)\rightarrow P(Y)$ be operations on $e^\star O(Y,\ \sigma)$ .

\defn\em A function $f$ : $(X,\ \tau)\rightarrow(Y,\ \sigma)$ is said to be $(e^\star_{[\gamma,\gamma']}, e^\star_{[\beta,\beta']})$-continuous if for each point $x\in X$ and each $e^\star$-open sets $W$ and $S$ of $(Y,\ \sigma)$ containing $f(x)$ there exist $e^\star$-open sets $U$ and $V$ of $(X,\ \tau)$ containing $x$ such that $f(U^{\gamma}\cap V^{\gamma})\subseteq W^{\beta}\cap S^{\beta'}.$

\thm\em\label{t4.1} Let (1), (2), (3), (4), (5), (6), and (7) be the following properties for a function $f$ : $(X, \tau) \rightarrow (Y, \sigma)$.
\begin{enumerate}[(1)]
	\item $f$ : $(X, \tau) \rightarrow (Y, \sigma)$ is $(e^\star_{[\gamma, \gamma']}, e^\star_{[\beta, \beta']})$-continuous.
	\item $f(e^\star Cl_{[\gamma, \gamma']}(A)) \subseteq e^\star Cl_{[\beta, \beta']}(f(A))$ for every subset $A$ of $(X, \tau)$.
	\item $e^\star Cl_{[\gamma, \gamma']}(f^{-1}(B)) \subseteq f^{-1}(e^\star Cl_{[\beta, \beta']}(B))$ for every subset $B$ of $(Y, \sigma)$.
	\item $f^{-1}(B)$ is $e^\star_{[\gamma, \gamma']}$-closed for every $e^\star_{[\beta, \beta']}$-closed set $B$ of $(Y, \sigma)$.
	\item $f(e^\star_{[\gamma, \gamma']}$-$Cl(A)) \subseteq e^\star_{[\beta, \beta']}$-$Cl(f(A))$ for every subset $A$ of $(X, \tau)$.
	\item $f^{-1}(V)$ is $e^\star_{[\gamma, \gamma']}$-open for every $e^\star_{[\beta, \beta']}$-open set $V$ of $(Y, \sigma)$.
	\item For each point $x \in X$ and each $e^\star_{[\beta, \beta']}$-open $W$ of $(Y, \sigma)$ containing $f(x)$, there exists $e^\star_{[\gamma, \gamma']}$-open $U$ of $(X, \tau)$ containing $x$ such that $f(U) \subseteq W.$
\end{enumerate}

Then (1) $\Rightarrow(2)\Leftrightarrow(3)\Rightarrow(4)\Leftrightarrow(5)\Leftrightarrow(6)\Leftrightarrow(7)$ hold.

\textbf{Proof.} (1) $\Rightarrow$ (2). Let $x\in e^\star Cl_{[\gamma,\gamma']}(A)$ and $W, S$ be $e^\star$-open sets of $(Y,\ \sigma)$ containing $f(x)$ . There exist $e^\star$-open sets $U$ and $V$ of $(X,\ \tau)$ containing $x$ such that $f(U^{\gamma}\cap V^{\gamma'})\subseteq W^{\beta}\cap S^{\beta'}$. Since $x\in e^\star Cl_{[\gamma,\gamma']}(A)$ , then $(U^{\gamma}\cap V^{\gamma'})\cap A\neq\phi$, implies that $ f(U^{\gamma}\cap V^{\gamma'})\cap f(A)\neq\phi$. Therefore, we have $ f(A)\cap(W^{\beta}\cap S^{\beta'})\neq\phi$. Therefore $f(x)\in e^\star Cl_{[\beta,\beta']}(f(A))$ , which implies that $x\in f^{-1}(e^\star Cl_{[\beta,\beta']}(f(A)))$ . Hence $e^\star Cl_{[\gamma,\gamma']}(A)\subseteq f^{-1}(e^\star Cl_{[\beta,\beta']}(f(A)))$ , so that $f(e^\star Cl_{[\gamma,\gamma']}(A))\subseteq e^\star Cl_{[\beta,\beta']}(f(A))$ .

(2) $\Rightarrow(3)$ . Let $B$ be any subset of $Y$. Then $f^{-1}(B)$ is a subset of $X$. By (2), we have $f(e^\star Cl_{[\gamma,\gamma']}(f^{-1}(B)))\subseteq e^\star Cl_{[\beta,\beta']}(f(f^{-1}(B)))\subseteq e^\star Cl_{[\beta,\beta']}(B)$ . Hence $e^\star Cl_{[\gamma,\gamma']}(f^{-1}(B))\subseteq f^{-1}(e^\star Cl_{[\beta,\beta']}(B))$ .

(3) $\Rightarrow$ (2). Let $A$ be any subset of $X$. Then $f(A)$ is a subset of $Y$. By (3), we have $e^\star Cl_{[\gamma,\gamma']}(f^{-1}f((A)))\subseteq f^{-1}(e^\star Cl_{[\beta,\beta']}(f(A)))$ . This implies that $e^\star Cl_{[\gamma,\gamma']}(A)\subseteq f^{-1}(e^\star Cl_{[\beta,\beta']}(f(A)))$ . Hence $f(e^\star Cl_{[\gamma,\gamma']}(A))\subseteq e^\star Cl_{[\beta,\beta']}(f(A))$ .

(3) $\Rightarrow(4)$ . Let $B$ be an $e^\star_{[\beta,\beta']}$-closed set of $(Y,\ \sigma)$ . By (3) and Theorem \ref{t3.1}, $e^\star Cl_{[\gamma,\gamma']}(f^{-1}(B))\subseteq f^{-1}(B)$ and hence $f^{-1}(B)$ is $e^\star_{[\gamma,\gamma']}$-closed.

(4) $\Rightarrow(5)$ . Let $A$ be any subset of $X$. Then $f(A) \subseteq e^\star_{[\beta, \beta']}$-$Cl(f(A))$ and $e^\star_{[\beta, \beta']}$-$Cl(f(A))$ is an $e^\star_{[\beta, \beta']}$-closed set in $Y$. Hence $A \subseteq f^{-1}(e^\star_{[\beta, \beta']}$-$Cl(f(A))).$ By (4), we have $f^{-1}(e^\star_{[\beta, \beta']}$-$Cl(f(A)))$ which is an $e^\star_{[\gamma, \gamma']}$-closed set in $X$. Therefore, $e^\star_{[\gamma, \gamma']}$-$Cl(A) \subseteq f^{-1}(e^\star_{[\beta, \beta']}$-$Cl(f(A))).$ Hence $f(e^\star_{[\gamma, \gamma']}$-$Cl(A)) \subseteq e^\star_{[\beta, \beta']}$-$Cl(f(A)).$

(5) $\Rightarrow(4)$ .Let $B$ be an $e^\star_{[\beta, \beta']}$-closed set of $(Y, \sigma)$. By (5), $e^\star_{[\gamma, \gamma']}$-$Cl(f^{-1}(B)) \subseteq f^{-1}(e^\star_{[\beta, \beta']}$-$Cl(f(f^{-1}(B)))) \subseteq f^{-1}(e^\star_{[\beta, \beta']}$-$Cl(B)) \subseteq f^{-1}(B)$. Therefore, by Proposition \ref{p3.13}, $f^{-1}(B)$ is $e^\star_{[\gamma, \gamma']}$-closed.

(5) $\Leftrightarrow(6)$ . This follows from Definition \ref{d34} and the equivalence of (4) $\Leftrightarrow(5)$ .

(6) $\Rightarrow$ (7). Let $W$ be any $e^\star_{[\beta, \beta']}$-open set in $Y$ containing $f(x)$, so its inverse image is an $e^\star_{[\gamma, \gamma']}$-open set in $X$. Since $f(x) \in W$, then $x \in f^{-1}(W)$ and by hypothesis $f^{-1}(W)$ is an $e^\star_{[\gamma, \gamma']}$-open set in $X$ containing $x$ so that $f(f^{-1}(W)) \subseteq W.$

(7) $\Rightarrow(6)$ . Let $V\in e^\star O(Y,\ \sigma)_{[\beta,\beta']}$. For each $x\in f^{-1}(V)$ , by (7), there exists an $e^\star_{[\gamma,\gamma']}$-open set $U_{x}$ containing $x$ such that $f(U_{x})\subseteq V$. Then we have $f^{-1}(V)=\bigcup\{U_{x}\in e^\star O(X,\ \tau)_{[\gamma,\gamma']}\ :\ x\in f^{-1}(V)\}$ and hence $f^{-1}(V)\in e^\star O(X,\ \tau)_{[\gamma,\gamma']}$ using Proposition \ref{p3.1} $\square $

\cor\em\label{c4.1} If $(Y, \sigma)$ is an $e^\star_{[\beta, \beta']}$-regular space, or operations $\beta$ and $\beta'$ are $e^\star$-open on $e^\star O(Y, \sigma)$, then all properties of Theorem \ref{t4.1} are equivalent.

\textbf{Proof.} By Theorem \ref{t4.1}, it is sufficient to prove the implication (4) $\Rightarrow(1)$ , where (1) and (4) are the properties of Theorem \ref{t4.1}

First, we show the implication under the assumption that $(Y,\ \sigma)$ is an $e^\star_{[\beta,\beta']^{-}}$ regular space. Let $x\in X$ and $W, S$ be $e^\star$-open sets of $(Y,\ \sigma)$ containing $f(x)$ . By Proposition \ref{p3.10}, $Y\backslash (W\cap S)$ is $e^\star_{[\beta,\beta']}$-closed. Then, $f^{-1}(Y\backslash (W\cap S))$ is $e^\star_{[\gamma,\gamma']}$-closed by (4) and hence $f^{-1}(W\cap S)$ is $e^\star_{[\gamma,\gamma']}$-open set of $(X,\ \tau)$ containing $x$. Therefore, there exist $e^\star$-open sets $U$ and $V$ of $(X,\ \tau)$ containing $x$ such that $U^{\gamma}\cap V^{\gamma'}\subseteq f^{-1}(W\cap S)$ and so $f(U^{\gamma}\cap V^{\gamma'})\subseteq W^{\beta}\cap S^{\beta'}$. This implies that $f$ is $(e^\star_{[\gamma,\gamma']},\ e^\star_{[\beta,\beta']})$-continuous.

Second, we suppose that the operations $\beta$ and $\beta'$ are $e^\star$-open. Let $x\in X$ and $W, S$ be $e^\star$-open sets of $(Y,\ \sigma)$ containing $f(x)$ . By using $e^\star$-openness of $\beta$ and $\beta'$, there exist an $e^\star_{\beta}$-open set $A$ and an $e^\star_{\beta'}$-open set $B$ such that $f(x)\in A\cap B$ and $A\cap B\subseteq W^{\beta}\cap S^{\beta'}$. By Proposition \ref{p3.6}, $Y\backslash (A\cap B)$ is $e^\star_{[\gamma,\gamma']}$-closed and hence $f^{-1}(Y\backslash (A\cap B))$ is $e^\star_{[\beta,\beta']}$-closed. Therefore, there exist $e^\star$-open sets $U$ and $V$ of $(X,\ \tau)$ containing $x$ such that $U^{\gamma}\cap V^{\gamma'}\subseteq f^{-1}(A\cap B)$ and so $f(U^{\gamma}\cap V^{\gamma'})\subseteq W^{\beta}\cap S^{\beta'}$. This implies that $f$ is $(e^\star_{[\gamma,\gamma']},\ e^\star_{[\beta,\beta']})$-continuous. $\square $

\rmk\em It is clear that the identity operation and the operation $X$ : $e^\star O(X,\ \tau)\rightarrow(X,\ \tau)$ are $e^\star$-open on $e^\star O(X,\ \tau)$ . Therefore, by Corollary \ref{c4.1}, if $\beta$ and $\beta'$ are chosen from mentioned operations above, then all properties of Theorem \ref{t4.1} are equivalent.

\rmk\em The converse of implication (1) $\Rightarrow(4)$ in Theorem \ref{t4.1} is not true in general as shown by the following example.

\ex\em\label{e4.1} Let $X=\{s_1,s_2,s_3\}$ and $\tau=\{\phi,\ X,\ \{s_1\},\ \{s_2\},\ \{s_1,s_2\}\}$ be a topology on $X$. Let $f$ : $(X,\ \tau)\rightarrow(X,\ \tau)$ be the identity. Let $\gamma=\gamma'=\beta'=X$: $e^\star O(X,\ \tau)\rightarrow P(X)$ be the operations on $e^\star O(X,\ \tau)$ and $\beta$ the closure operation on $e^\star O(X,\ \tau)$ . Then, the condition (4) in Theorem \ref{t4.1} is true. It is shown that $f$ is not $(e^\star_{[\gamma,\gamma']},\ e^\star_{[\beta,\beta']})$-continuous.

Theorem \ref{t4.1} suggests the following.

\rmk\em\label{r4.3} If $f:(X,\ \tau)\rightarrow(Y,\ \sigma)$ is $(e^\star_{[\gamma,\gamma']},\ e^\star_{[\beta,\beta']})$-continuous, then the induced function $f:(X,\ e^\star O(X,\ \tau)_{[\gamma,\gamma']})\rightarrow(Y,\ e^\star O(Y,\ \sigma)_{[\beta,\beta']})$ is continuous.

\rmk\em The converse of Remark \ref{r4.3} is not true in general as shown by the following example.

\ex\em Let $(X,\ \tau)$ as in Example \ref{e4.1} and $f:(X,\ \tau)\rightarrow(X,\ \tau)$ be a function defined by

$f(x)=\left\{\begin{array}{l}
s_2\ \mathrm{i}\mathrm{f}\ x=s_1,\\
s_3\ \mathrm{i}\mathrm{f}\ x=s_2,\\
s_1\ \mathrm{i}\mathrm{f}\ x=s_3,
\end{array}\right.$

moreover let $\gamma=\beta$ be the closure operation on $e^\star O(X,\ \tau)$ and $\gamma'=\beta'= X$ : $e^\star O(X,\ \tau)\rightarrow P(X)$ . Then, $e^\star O(X,\ \tau)_{[\gamma,X]} =\{\phi,\ X\}$ and it is shown that $f$ is not $(e^\star_{[\gamma,X]},\ e^\star_{[\gamma,X]})$-continuous. However, $f$ : $(X,\ e^\star O(X,\ \tau)_{[\gamma,\gamma']})\rightarrow (Y,\ e^\star O(Y,\ \sigma)_{[\beta,\beta']})$ is continuous.

Let $(X,\ \tau), (Y,\ \sigma)$ and $(Z,\ \eta)$ be spaces and $\gamma, \gamma'$ : $e^\star O(X,\ \tau)\rightarrow P(X)$ , $\beta, \beta'$ : $e^\star O(Y)\rightarrow P(Y)$ and $\delta, \delta'$ : $e^\star O(Z)\rightarrow P(Z)$ , be operations on $e^\star O(X,\ \tau)$ , $e^\star O(Y,\ \sigma)$ and $e^\star O(Z,\ \eta)$ , respectively.

\thm\em If $f$ : $(X, \tau) \rightarrow (Y, \sigma)$ is $(e^\star_{[\gamma, \gamma']}, e^\star_{[\beta, \beta']})$-continuous and $g$ : $(Y, \sigma) \rightarrow (Z, \eta)$ is $(e^\star_{[\beta, \beta']}, e^\star_{[\delta, \delta']})$-continuous, then its composition $g\circ f$ : $(X, \tau) \rightarrow (Z, \eta)$ is $(e^\star_{[\gamma, \gamma']}, e^\star_{[\delta, \delta']})$-continuous.

\textbf{Proof.} Let $x\in X, K$ and $L$ be $e^\star$-open sets of $Z$ containing $g(f(x))$ . Since $g$ is $(e^\star_{[\beta,\beta']},\ e^\star_{[\delta,\delta']})$-continuous, then there exist $e^\star$-open sets $W$ and $S$ of $Y$ containing $f(x)$ such that $g(W^{\beta}\cap S^{\beta'})\subseteq K^{\delta}\cap L^{\delta'}$. Also, since $f$ is $(e^\star_{[\gamma,\gamma']}, e^\star_{[\beta,\beta']})$-continuous, then there exist $e^\star$-open sets $U$ and $V$ of $X$ containing $x$ such that $f(U^{\gamma}\cap V^{\gamma'})\subseteq W^{\beta}\cap S^{\beta'}$. This implies that $ f(U^{\gamma}\cap V^{\gamma'})\subseteq W^{\beta}\cap S^{\beta'}\subseteq g^{-1}(K^{\delta}\cap L^{\delta'})$ . Then we obtain $(g\circ f)(U^{\gamma}\cap V^{\gamma'})\subseteq K^{\delta}\cap L^{\delta'}$. Therefore, $g\circ f$ is $(e^\star_{[\gamma,\gamma']},\ e^\star_{[\delta,\delta']})$-continuous. $\square $

\defn\em A function $f$ : $(X,\ \tau)\rightarrow(Y,\ \sigma)$ is said to be $(e^\star_{[\gamma,\gamma']}, e^\star_{[\beta,\beta']})$-closed if for $e^\star_{[\gamma,\gamma']}$-closed set $A$ of $X, f(A)$ is $e^\star_{[\beta,\beta']}$-closed in $Y.$

\prop\em Let $f$ : $(X, \tau) \rightarrow (Y, \sigma)$ be an $(e^\star_{[\gamma, \gamma']}, e^\star_{[\beta, \beta']})$-closed function. Then, for each subset $B$ of $(Y, \sigma)$ and each $e^\star_{[\gamma, \gamma']}$-open set $U$ containing $f^{-1}(B)$, there exists an $e^\star_{[\beta, \beta']}$-open set $V$ such that $B \subseteq V$ and $f^{-1}(V) \subseteq U.$

\textbf{Proof.} Let $V=Y\backslash f(X\backslash U)$ . Then $V$ is $e^\star_{[\beta,\beta']}$-open. Thus $f^{-1}(B)\subseteq U$ implies $B\subseteq V$ and $ f^{-1}(V)=f^{-1}(Y\backslash f(X\backslash U))=X\backslash f^{-1}(f(X\backslash U))\subseteq X\backslash (X\backslash U)=U$, or $f^{-1}(V)\subseteq U. \square $

\prop\em If $f:(X, \tau) \rightarrow (Y, \sigma)$ is bijective and $ f^{-1}:(Y, \sigma) \rightarrow (X, \tau)$ is $(e^\star_{[\beta, \beta']}, e^\star_{[\gamma, \gamma']})$-continuous, then $f$ is $(e^\star_{[\gamma, \gamma']}, e^\star_{[\beta, \beta']})$-closed.

\textbf{Proof.} This follows from definitions and Theorem \ref{t4.1}

\end{document}